\documentclass[11pt, reqno]{amsart}
\usepackage{amssymb}

\newtheorem{theorem}{Theorem}[section]

\newtheorem{proposition}[theorem]{Proposition}
\newtheorem{lemma}[theorem]{Lemma}
\newtheorem{corollary}[theorem]{Corollary}

\theoremstyle{definition} 
\newtheorem{definition}[theorem]{Definition}
\newtheorem{example}[theorem]{Example}
\newtheorem{remark}[theorem]{Remark}
\newtheorem{question}[theorem]{Question}

\newcommand{\R}{\mathbb R}
\newcommand{\T}{\mathbb T}
\newcommand{\Z}{\mathbb Z}
\newcommand{\kk}{\mathbb k}

\DeclareMathOperator{\sys}{{\rm sys}}

\DeclareMathOperator{\vol}{{\rm vol}}

\DeclareMathOperator{\length}{{\rm length}}

\DeclareMathOperator{\card}{{\rm card}}

\DeclareMathOperator{\MinEnt}{{\rm MinEnt}}
\DeclareMathOperator{\Aut}{{\rm Aut}}
\DeclareMathOperator{\Ent}{{\rm Ent}}
\DeclareMathOperator{\Hom}{{\rm Hom}}
\DeclareMathOperator{\Diam}{{\rm Diam}}

\numberwithin{equation}{section} 
\numberwithin{figure}{section}

\begin{document}

\author[S.~Sabourau]{St\'ephane Sabourau} \address{Laboratoire de
Math\'ematiques et Physique Th\'eorique, Universit\'e de Tours, Parc de Grandmont, 37400
Tours, France} \email{sabourau@lmpt.univ-tours.fr}

\thanks{\emph{Current address:} 
Department of Mathematics,
University of Pennsylvania,
209 South 33rd Street,
Philadelphia, PA 19104-6395, USA}
\thanks{\emph{E-mail:} \tt{sabourau@lmpt.univ-tours.fr}}




\title[Systolic volume and minimal entropy]{Systolic volume and minimal entropy of aspherical manifolds}

\begin{abstract}
We establish isosystolic inequalities for a class of manifolds which includes the aspherical manifolds. In particular, we relate the systolic volume of aspherical manifolds first to their minimal entropy, then to the algebraic entropy of their fundamental groups. 
\end{abstract}

\maketitle

\tableofcontents

\section{Introduction}

Let $M$ be a nonsimply connected closed $n$-manifold endowed with a Riemannian metric~$g$. The (homotopy) systole of~$(M,g)$, denoted $\sys(M,g)$, is defined as the length of the shortest noncontractible loop in~$M$.
The optimal systolic volume of $M$, denoted $\sigma(M)$, is defined as
\begin{equation} \label{eq:sigma}
\sigma(M) = \inf_g \frac{\vol(M,g)}{\sys(M,g)^n}
\end{equation}
where $g$ runs over the space of all metrics on~$M$.

The study of the systolic volume constitutes the heart of systolic geometry. We refer the reader to the expository texts~\cite{ber}, \cite{gro96}, \cite{gro99} and~\cite{c-k} and the references therein for an account on the subject.

In 1983, M.~Gromov proved the following fundamental systolic inequality whose converse for orientable manifolds was established by I.~Babenko in~\cite{bab}.

\begin{theorem}[\cite{gro83}] \label{theo:ess}
There exists a positive constant $C_n$ such that every $n$-dimensional essential manifold~$M$ satisfies
\begin{equation} \label{eq:ess}
\sigma(M) \geq C_n.
\end{equation}
\end{theorem}

Recall that a closed manifold~$M$ is said to be essential if there is a map~$f$ from~$M$ into a~$K(\pi,1)$ such that $f_*([M]_\kk) \neq 0$ in~$H_n(\pi;\kk)$ where $\kk = \Z$ if~$M$ is orientable or $\Z_2$ otherwise. In particular, aspherical manifolds are essential.

This result extends earlier works on surfaces (see the references mentioned above).
Note that the precise value of the optimal systolic volume of essential manifolds is known only for the 2-torus (C.~Loewner, unpublished, \cite{gro99}), the projective plane (P.M.~Pu, \cite{pu}) and the Klein bottle (C.~Bavard, \cite{bav}, \cite{sak}).

The systolic inequality~\eqref{eq:ess} can be improved by taking into account the topology of the manifold.
For instance, closed surfaces of large genus have a large systolic volume (see~\cite[5.3]{gro83}).
More generally, M.~Gromov showed that the systolic volume of a closed manifold~$M$ is bounded from below in terms of its simplicial volume~$||M||$, its simplicial height~$h(M)$ and, in the aspherical case, the sum of its Betti numbers $b(M) = \sum_{i} b_{i}(M)$ (see~\cite[3.C]{gro96} for precise definitions).
Specifically, we have

\begin{theorem}[{\cite[6.4.C''' \& D']{gro83}, \cite[3.C]{gro96} and also \cite[4.46]{gro99}}] \label{theo:simpl}
Every $n$-dimensional closed manifold~$M$ satisfies
\begin{equation} \label{eq:grosimpl}
\sigma(M) \geq C_n \frac{||M||}{\log^n(1+||M||)},
\end{equation}
\begin{equation} \label{eq:groheight}
\sigma(M) \geq C_n \frac{h(M)}{\exp(C'_{n} \sqrt{\log h(M)})}.
\end{equation}
Furthermore, if $M$ is aspherical, we have
\begin{equation} \label{eq:grobetti}
\sigma(M) \geq C_n \frac{b(M)}{\exp(C'_{n} \sqrt{\log b(M)})}.
\end{equation}
Here, the positive constants $C_{n}$ and~$C'_{n}$ depend only on the dimension of~$M$.
\end{theorem}

Note that inequality~\eqref{eq:grobetti} is a consequence of~\eqref{eq:groheight} since \mbox{$h(M) \leq b(M)$} for aspherical manifolds.

In \cite{bab}, \cite{bab04a} and~\cite{bab04}, I.~Babenko compared the optimal systolic volumes of the manifolds with the same fundamental group. Namely, he proved that the optimal systolic volume is a homotopy invariant (see~\cite{bab}). He also showed that in a large number of cases, the optimal systolic volume of a manifold depends only on its homology class in the corresponding Eilenberg-MacLane space (see~\cite{bab04a} and~\cite{bab04} for a more precise statement). \\

In this article, we prove two other systolic inequalities for aspherical manifolds using, on the one hand, the minimal entropy and, on the other hand, the algebraic entropy of the fundamental groups. Loosely speaking, we show that aspherical manifolds with ``complicated'' fundamental groups have large systolic volume. \\

In order to state our first result, we need to introduce some definitions.

Given a closed $n$-dimensional manifold~$(M,g)$, denote by~$(\tilde{M}, \tilde{g})$ the universal Riemannian covering .
Fix $x_0 \in M$ and a lift $\tilde{x}_0 \in \tilde{M}$ of~$x_0$.

The volume entropy (or asymptotic volume) of $(M,g)$ is defined as
\begin{equation} \label{eq:defent}
\Ent(M,g) = \lim_{R \rightarrow + \infty} \frac{\log(\vol_{\tilde{g}} B(\tilde{x}_0,R) )}{R}
\end{equation}
where $\vol_{\tilde{g}} B(\tilde{x}_0,R)$ is the volume of the ball centered at $\tilde{x}_0$ with radius~$R$ in~$\tilde{M}$. Since $M$ is compact, the limit in \eqref{eq:defent} exists and does not depend on the point $\tilde{x}_0 \in \tilde{M}$ (see~\cite{man}).
This asymptotic invariant describes the exponential growth rate of the volume on the universal covering.

Define the minimal entropy of~$M$ as
\begin{equation}
\MinEnt(M) = \inf_g \Ent(M,g) \vol(M,g)^{\frac{1}{n}}
\end{equation}
where $g$ runs over the space of all metrics on~$M$.

M.~Gromov showed in~\cite{gro81} that the simplicial volume is related to the minimal entropy. More precisely, there exists a positive constant~$C_n$ such that every closed $n$-manifold~$M$ with simplicial volume~$||M||$ satisfies
\begin{equation} \label{eq:gro}
\MinEnt(M)^n \geq C_n ||M||.
\end{equation}
When $M$ admits a negatively curved locally symmetric metric~$g_0$, it has been proved by A.~Katok in~\cite{kat} (for $n=2$) and by G.~Besson, G.~Courtois and S.~Gallot in~\cite{bcg95} (for $n \geq 3$) that
\begin{equation} \label{eq:bcg}
\MinEnt(M)^n  = \Ent(g_0)^n \vol(g_0).
\end{equation}

I.~Babenko showed in~\cite{bab} that both the optimal systolic volume and the minimal volume entropy are homotopy invariants of manifolds.

The first result of this article shows how these two invariants are related.
More precisely, we prove the following partial generalization of inequality~\eqref{eq:grosimpl} in Theorem~\ref{theo:simpl}.

\begin{theorem} \label{theo:ent}
Let $\Phi:M \longrightarrow K$ be a degree~1 map between two oriented closed manifolds of the same dimension~$n \neq 3$ such that $K$ is a $K(\pi_1(M),1)$-space. Then, there exists a positive constant $C_n$ depending only on~$n$ such that
\begin{equation} \label{eq:ent}
\sigma(M) \geq C_n \frac{\MinEnt(M)^n}{\log^n(1+\MinEnt(M))}.
\end{equation}
\end{theorem}

If~$M$ is aspherical, we can take the identity map for~$\Phi$. This immediately leads to the following result.

\begin{corollary} \label{cor:ent}
Every orientable aspherical closed manifold of dimension~$n$ different from~$3$ satisfies the systolic inequality~\eqref{eq:ent}.
\end{corollary}

Theorem~\ref{theo:ent} is more general than Corollary~\ref{cor:ent}.
Indeed, the connected sum of an aspherical manifold satisfying the hypothesis of Theorem~\ref{theo:ent} with a simply connected orientable closed manifold nonhomeomorphic to the sphere again satisfies the hypothesis of Theorem~\ref{theo:ent} but is not aspherical.

Inequality~\eqref{eq:ent} combined with~\eqref{eq:gro} implies inequality~\eqref{eq:grosimpl}.
Thus, Theorem~\ref{theo:ent} sharpens Theorem~\ref{theo:simpl}~\eqref{eq:grosimpl} for aspherical manifolds of dimension different from~$3$, but it does not apply to arbitrary closed manifolds as Theorem~\ref{theo:simpl}~\eqref{eq:grosimpl} does.
It is an open question whether inequality~\eqref{eq:ent} holds for every closed manifold.

The restriction on the dimension of~$M$ in Theorem~\ref{theo:ent} is due to technical reasons: in the proof, we use a construction of I.~Babenko~\cite{bab04} which relies on obstruction theory and requires that $n$ is different from~$3$ (see Section~\ref{sec:minent} for further explanation). \\

Let us now introduce the algebraic entropy of a finitely generated group.

Let $\Sigma$ be a finite generating set of a group~$\Gamma$. 
The algebraic (or word) length of an element $\alpha \in \Gamma$ with respect to~$\Sigma$ is noted~$| \alpha |_\Sigma$. It is defined as the smallest integer $k \geq 0$ such that $\alpha = \alpha_1 \dots \alpha_k$ where $\alpha_i \in \Sigma \cup \Sigma^{-1}$. By definition, the neutral element~$e$ is the only element of $\Gamma$ with null algebraic length.

The algebraic entropy of $\Gamma$ with respect to~$\Sigma$ is defined as
\begin{equation} \label{eq:entalg}
\Ent_{alg}(\Gamma,\Sigma) = \limsup_{R \rightarrow + \infty} \frac{\log(N_\Sigma(R))}{R}.
\end{equation}
where $N_\Sigma(R)= \card \{ \alpha \in \Gamma \mid |\alpha|_\Sigma \leq R \}$ is the cardinal of the $R$-ball of~$(\Gamma,|.|_{\Sigma})$ centered at its origin.
If $\Gamma$ is the fundamental group of a closed polyhedron, then the sup-limit is a limit (see~\cite{man}).

Define the minimal algebraic entropy of~$\Gamma$ as
\begin{equation}
\Ent_{alg}(\Gamma) = \inf_\Sigma \Ent_{alg}(\Gamma,\Sigma)
\end{equation}
where $\Sigma$ runs over the space of all generating sets of~$\Gamma$.

The second result of this article shows how the systolic volume of some manifolds is related to the algebraic entropy of their fundamental groups.
More precisely, we have

\begin{theorem} \label{theo:alg}
Let $\Phi:M \longrightarrow K$ be a degree~1 map between two $n$-dimensional oriented closed manifolds such that $K$ is a $K(\pi_1(M),1)$-space.
Then, there exists a positive constant $C_n$ depending only on~$n$ such that
\begin{equation} \label{eq:alg}
\sigma(M) \geq C_n \frac{\Ent_{alg}(\pi_1(M))}{\log(1+\Ent_{alg}(\pi_1(M)))}.
\end{equation}
\end{theorem}
Note that the dimension~$n$ can be equal to~3 in this result.
As previously, we have the following particular result.

\begin{corollary}
Every $n$-dimensional orientable aspherical closed manifold satisfies the systolic inequality~\eqref{eq:alg}.
\end{corollary}

As shown in Example~\ref{ex:handles}, inequality~\eqref{eq:alg} does not hold for every orientable closed manifold: extra topological conditions are required, as in Theorem~\ref{theo:alg}.

The constants~$C_n$ we obtain in inequalities~\eqref{eq:ent} and \eqref{eq:alg} can be explicitly computed.
Theorems~\ref{theo:ent} and~\ref{theo:alg} are consequences of more general results - namely, Theorems~\ref{theo:Gsys} and~\ref{theo:Galg} - stated in terms of relative systole and relative entropy (see Section~\ref{sec:rel} for a definition of these notions).
Note that their proofs are not based on filling invariant estimates. \\

Let us now mention some other related works, though of a quite different nature, on universal systolic inequalities.

Upper bounds on the optimal systolic volume of the connected sums and the cyclic coverings of manifolds can be found in~\cite{b-s} for surfaces and in~\cite{b-b} in higher dimensions.

The (homotopy) systole can be replaced by stable or conformal systoles of any dimension.
In this case, multiplicative relations in the cohomology ring lead to systolic inequalities (see~\cite{gro83}, \cite{gro99}, \cite{b-k1}, \cite{b-k2}, \cite{k03}, \cite{i-k}, \cite{bcik}, \cite{k-r} for more precise statements and related results).
Note that some of these stable/conformal systolic inequalities are sharp.

In contrast with these results, unstable higher dimensional systoles typically satisfy no nontrivial systolic inequalities (see~\cite{f99}, \cite{ks01} and~\cite{bab02} for the most general and recent results). \\

This article is organized as follows.
The definition of pseudomanifolds and the notions of relative systole and relative entropy are presented in Section~\ref{sec:rel}.
In Section~\ref{sec:minent}, we show that under some topological conditions, the relative minimal entropy of geometric cycles representing a manifold is not less than the relative minimal entropy of this manifold. The proof of this comparison result makes use of surgery and extension procedures and relies on I.~Babenko's constructions.
In Section~\ref{sec:ent}, we introduce regular geometric cycles whose existence and main properties have been established by M.~Gromov. Then, we bound from above their relative minimal entropy in terms of their relative systolic volume.
Combining the results of Sections~\ref{sec:minent} and~\ref{sec:ent}, we derive in Section~\ref{sec:Gsys} a lower bound on the relative systolic volume of some manifolds in terms of their relative minimal entropy.
In Section~\ref{sec:Galg}, we show that under some topological conditions, the relative systolic volume of manifolds can be bounded from below in terms of the algebraic entropy of their fundamental group. 

\section{Pseudomanifolds, relative systole and relative entropy} \label{sec:rel}

In this section, we recall the definition of pseudomanifolds and introduce the notions of relative systole and relative entropy. \\

The definition of manifolds appears too restrictive in the study of the systolic volume.
Instead, we will consider Riemannian polyhedra, i.e., polyhedra endowed with a piecewise linear Riemannian metric (see~\mbox{\cite[\S2]{bab}} for a more precise definition).
Unless stated otherwise, all polyhedra will be finite.
By definition, a finite $n$-dimensional polyhedron $X$ is said to be orientable if $H_n(X;\Z) \simeq \Z$.

The class of closed orientable $n$-dimensional pseudomanifolds plays an important role in our study, especially in relation to inequalities~\eqref{eq:reg0} and~\eqref{eq:reg}.

\begin{definition}[\cite{spa}]
A closed $n$-dimensional pseudomanifold is a finite simplicial complex~$X$ such that
\begin{enumerate}
\item every simplex of~$X$ is a face of some $n$-simplex of~$X$;
\item every $(n-1)$-simplex of~$X$ is the face of exactly two $n$-simplices of~$X$;
\item given two $n$-simplices $s$ and $s'$ of~$X$, there exists a finite sequence $s=s_1, s_2, \dots, s_m=s'$ of $n$-simplices of~$X$ such that $s_i$ and $s_{i+1}$ have an $(n-1)$-face in common.
\end{enumerate}
\end{definition}

The $n^{\text{th}}$ homology group~$H_n(X;\Z)$ of a closed $n$-dimensional pseudomanifold is either isomorphic to~$\Z$ or trivial (see~\cite{spa}).

\begin{lemma} \label{lem:cell}
Every closed $n$-dimensional pseudomanifold admits a cellular decomposition with exactly one $n$-cell.
\end{lemma}

\begin{proof}
Consider a closed $n$-pseudomanifold~$X$ with its simplicial structure. There exists a tree~$T$ in~$X$ whose vertices agree with the barycenters of the $n$-simplices of~$X$ and edges are segments joining the barycenters of adjacent $n$-simplices. Here, we say that two $n$-simplices are adjacent if they have an $(n-1)$-face in common.
The $(n-1)$-simplices of~$X$ which do not intersect~$T$ form an $(n-1)$-dimensional complex~$X'$. It is possible to define a cellular complex homeomorphic to~$X$ by gluing an $n$-cell to~$X'$.
\end{proof}

Let us now extend the classical definitions of the systole and the entropy.

Let $X$ be a finite $n$-dimensional polyhedron endowed with a piecewise Riemannian metric~$g$. Let $\psi: \pi_1(X) \longrightarrow \pi$ be a group homomorphism with kernel $H \lhd \pi_1(X)$. Denote by $\psi_\sharp$ the map induced by~$\psi$ between the conjugacy classes.

The $\psi$-systole of $(X,g)$, denoted $\sys_\psi(X,g)$, is defined as the length of the shortest loop of~$X$ whose image by~$\psi_\sharp$ is nontrivial.

The $\psi$-systolic volume of $(X,g)$ is defined as
\begin{equation}
\sigma_\psi(X,g) = \frac{\vol(X,g)}{\sys_\psi(X,g)^n}
\end{equation}
and the optimal $\psi$-systolic volume of $X$ is defined as
\begin{equation} \label{eq:defsigpsi}
\sigma_\psi(X) = \inf_g \sigma_\psi(X,g)
\end{equation}
where $g$ runs over the space of all piecewise linear Riemannian metrics on~$X$.
Here, the volume of an $n$-dimensional Riemannian polyhedron~$(X,g)$ is defined as the sum of the $n$-volumes of all the $n$-simplices of~$X$.
Thus, the $n$-volume agrees with the $n$-dimensional Hausdorff measure.
Note that the optimal $\psi$-systolic volume of $X$ depends rather on the kernel~$H$ of~$\psi$ than on the homomorphism~$\psi$ itself. We clearly have
\begin{equation}
\sigma_\psi(X) \leq \sigma(X).
\end{equation}

Let $p: X_H \longrightarrow X$ be the covering corresponding to the normal subgroup~$H \lhd \pi_1(X)$.
By definition, the deck transformation group of~$X_H$, noted $\Aut(X_H)$, agrees with $\pi_1(X)/H$ and $\pi_1(X_H)$ is isomorphic to~$H$.
The lift of~$g$ to~$X_H$ is noted~$\overline{g}$.
Fix $x_0$ in $X$ and a lift $\overline{x}_0$ in $X_H$.

The volume entropy of $X$ relative to $H$ (or $H$-entropy) is defined as
\begin{equation} \label{eq:defentphi}
\Ent_H(X,g) = \lim_{R \rightarrow +\infty} \frac{1}{R} \log(\vol_{\overline{g}} B(\overline{x}_0,R))
\end{equation}
where $\vol_{\overline{g}} B(\overline{x}_0,R)$ is the volume of the $\overline{g}$-ball centered at~$\overline{x}_0$ of radius~$R$ in~$X_H$.
Since $X$ is compact, the limit in \eqref{eq:defentphi} exists and does not depend on the point $\overline{x}_0 \in X_H$ (see~\cite{man}).
When $H$ is trivial, the \mbox{$H$-entropy} agrees with the classical volume entropy.
Furthermore, if \mbox{$H_1 \lhd H_2 \lhd \pi_1(M)$}, then $\Ent_{H_1}(X,g) \geq \Ent_{H_2}(X,g)$. In particular, one has $\Ent_H(X,g) \leq \Ent(X,g)$.

The minimal $H$-entropy of~$X$ is defined as
\begin{equation} \label{eq:defentH}
\MinEnt_H(X) = \inf_g \Ent_H(X,g) \vol(X,g)^{\frac{1}{n}}
\end{equation}
where $g$ runs over the space of all piecewise linear Riemannian metrics on~$X$.
When $H$ is trivial, the minimal $H$-entropy agrees with the classical minimal  entropy.

The $H$-entropy can also be defined by using the following result.
\begin{lemma} \label{lem:entH}
Let $X$ be a closed polyhedron endowed with a piecewise linear Riemannian metric~$g$ and $\psi:\pi_1(X) \longrightarrow~\pi$ be a group homomorphism with kernel~$H$. Then,
\begin{equation}
\Ent_H(X,g) = \lim_{T \rightarrow +\infty} \frac{\log(P_H(T))}{T}
\end{equation}
where $P_H(T)$ is the number of classes in $\Aut(X_H) \simeq \pi_1(X,x_0)/H$ which can be represented by loops of length at most~$T$ based at some fixed point~$x_0$.
\end{lemma}

\begin{proof}
The group $\Gamma:=\Aut(X_H) \simeq \pi_1(X,x_0)/H$ acts
on~$X_H$ by isometries.  The orbit of~$\overline{x}_0$ by~$\Gamma$
is noted~$\Gamma.\overline{x}_0$ where $\overline{x}_0$ is a lift of~$x_0$ in $X_H$.  Consider a fundamental domain~$\Delta$
for the action of~$\Gamma$ containing~$\overline{x}_0$. Denote by~$D$ the
diameter of~$\Delta$.  The cardinal of $\Gamma.\overline{x}_0 \cap
B(\overline{x}_0,R)$ is bounded from above by the number of translated
fundamental domains~$\gamma.\Delta$ contained in $B(\overline{x}_0,R+D)$
and bounded from below by the number of translated fundamental
domains~$\gamma.\Delta$ contained in $B(\overline{x}_0,R)$.  Therefore,
we have
\begin{equation}
\frac{\vol(B(\overline{x}_0,R))}{\vol(X,g)} \leq \card
(\Gamma.\overline{x}_0 \cap B(\overline{x}_0,R)) \leq
\frac{\vol(B(\overline{x}_0,R+D))}{\vol(X,g)}.
\end{equation}
Take the log of these terms, multiply them by~$\frac{1}{R}$ and let $R$ go to infinity.  The left-hand term and the right-hand term both yield $\Ent_H(X,g)$.
The result follows since $P_H'(R) = \card (\Gamma.\overline{x}_0 \cap
B(\overline{x}_0,R))$.
\end{proof}



\section[Minimal entropy of geometric cycles]{Minimal entropy of geometric cycles representing aspherical manifolds} \label{sec:minent}

In this section, we compare the relative minimal entropy of some manifolds with the relative minimal entropy of their representing geometric cycles. \\

Let $\pi$ be a discrete group and $K(\pi,1)$ be an Eilenberg-MacLane space.
Let $h$ be a homology class in $H_n(\pi;\Z):= H_n(K(\pi,1);\Z)$.
A geometric cycle representing the homology class~$h$ is a map $\Psi:X \longrightarrow K(\pi,1)$ from a closed oriented $n$-dimensional Riemannian pseudomanifold~$(X,g)$ such that $\Psi_*[X]=h$ where $[X]$ is the fundamental class of~$X$.
The map~$\Psi$ induces a homomorphism $\psi:\pi_1(X) \longrightarrow \pi$ between the fundamental groups, whose kernel is noted~$H$.

The minimal entropy of $h \in H_n(\pi;\Z)$ is defined as
\begin{equation}
\MinEnt(h) = \inf_X \MinEnt_H(X)
\end{equation}
where $X=(X,\Psi,g)$ runs over all the geometric cycles representing~$h$.

\begin{theorem} \label{theo:minent}
Let $\pi$ be a discrete group.
Let $\Phi:M \longrightarrow K$ be a degree~1 map between two oriented closed manifolds of the same dimension~$n \neq 3$ such that $K$ is a $K(\pi,1)$-space.
Then,
\begin{equation}
\MinEnt_G(M) = \MinEnt(\Phi_*[M])
\end{equation}
where $G$ is the kernel of the homomorphism $\phi:\pi_1(M) \longrightarrow \pi$ induced by~$\Phi$ and $[K]$ is the fundamental class of~$K$.
\end{theorem}

\begin{remark} \label{rem:minent}
Under the assumptions of Theorem~\ref{theo:minent}, we have \mbox{$\Phi_*[M]=[K]$}. We can apply again Theorem~\ref{theo:minent} to the identity map of~$K$. This yields
\begin{equation}
\MinEnt_G(M) = \MinEnt(K).
\end{equation}
\end{remark}
\begin{remark}
Every homomorphism $\phi:\pi_1(M) \longrightarrow \pi$ induces a map $\Phi:M \longrightarrow K$ (uniquely defined up to homotopy) whose homomorphism induced between the fundamental groups agrees with~$\phi$.
Thus, if $M$ is aspherical, we can take $\pi=\pi_1(M)$, $K=K(\pi_1(M),1)=M$ and $\Phi$ the map induced by the identity homomorphism on~$\pi_1(M)$. In this case, we obtain
\begin{equation}
\MinEnt(M) = \MinEnt([M]).
\end{equation}
\end{remark}
\begin{remark}
Under the hypotheses of Theorem~\ref{theo:minent}, the group~$\pi$ is necessarily a quotient of~$\pi_1(M)$ (see Lemma~\ref{lem:surj}).
\end{remark}

Let $X$ be a closed $n$-dimensional polyhedron and $\psi:\pi_1(X) \longrightarrow \pi$ be a group homomorphism with kernel~$H$.
An $H$-extension of $X$ is a closed $n$-dimensional polyhedron~$X'$ obtained by attaching to $X$ a finite number of $k$-cells with $2 \leq k \leq n-1$ along simplicial maps $h:S^{k-1} \longrightarrow X'^{(k-1)}$ into the $(k-1)$-skeleton of~$X'$, with $[h] \in H$ when $k=2$.
In particular, an $H$-extension $X'$ of~$X$ contains~$X$ and may be written as
\begin{equation} \label{eq:X'}
X' = X \bigcup_{k=2}^{n-1} \bigcup_{i_k =1}^{p_k} B_{i_k}^k.
\end{equation}
The group $\pi_1(X')$ is isomorphic to a quotient of~$\pi_1(X)$ by a normal subgroup of~$H$. Thus, the map~$\psi$ induces a homomorphism \mbox{$\psi':\pi_1(X') \longrightarrow \pi$} whose kernel~$H'$ is a quotient of~$H$.

The inclusion map $i:X \hookrightarrow X'$ induces an isomorphism between $H_n(X;\Z)$ and $H_n(X';\Z)$. In particular, if $X$ is an orientable polyhedron, then the same holds for~$X'$. In this case, the degree of the inclusion map $i:X \hookrightarrow X'$ equals~1 (i.e., $i_*[X]=[X']$ in $H_n(X';\Z)$ for corresponding orientations). Recall that the degree of a map $f$ between two oriented $n$-dimensional polyhedra $X$ and $X'$ is defined as the unique integer $p$ such that $f_*[X]=p[X']$ in $H_n(X';\Z)$.
When $H$ is trivial, the space~$X'$ is simply called an extension of~$X$.

A simplicial map $f:X \longrightarrow Y$ between two finite simplicial polyhedra is said to be $n$-monotone if the preimage of every open $n$-simplex of~$Y$ is either an open $n$-simplex of~$X$ or an empty set. \\

The following result is due to I.~Babenko (see~\cite[\S2]{bab}).

\begin{lemma} \label{lem:mon}
Let $X_i$, $i=1,2$, be two closed $n$-dimensional polyhedra and $\psi_i:\pi_1(X_i) \longrightarrow \pi$ be two group homomorphisms with kernel~$H_i$.
Assume that there exists an $n$-monotone map $f:X_1 \longrightarrow X_2$ such that $\psi_1 = \psi_2 \circ f_*$.
Then,
\begin{equation}
\MinEnt_{H_1}(X_1) \leq \MinEnt_{H_2}(X_2)
\end{equation}
\end{lemma}

\begin{proof}
Fix $\varepsilon >0$ and a metric~$g$ on~$X_2$.
From~\cite[\S2]{bab} (see also~\cite[Lemme~3.1]{bab04}), there exists a metric~$g_\varepsilon$ on~$X_1$ such that \mbox{$f:(X_1,g_\varepsilon) \longrightarrow (X_2,g)$} is nonexpanding and $\vol(X_1,g_\varepsilon) \leq \vol(X_2,g) + \varepsilon$.
Since $\psi_1 = \psi_2 \circ~f_*$, we have $f_*^{-1}(H_2) = H_1$.
Thus, $f_*$ induces an injective homomorphism between $\pi_1(X_1)/H_1$ and $\pi_1(X_2)/H_2$. Moreover, if $\gamma$ is a loop of~$X_1$ based at some fixed point~$x_0$, its image $f(\gamma)$ is a loop of~$X_2$ based at $f(x_0)$ such that $\length(f(\gamma)) \leq \length(\gamma)$.
Therefore, $P_{H_1}(T) \leq P_{H_2}(T)$.
This implies that $\Ent_{H_1}(X_1,g_\varepsilon) \leq \Ent_{H_2}(X_2,g)$ from Lemma~\ref{lem:entH}. The result follows.
\end{proof}

The following result is an immediate consequence of Lemma~\ref{lem:mon}.

\begin{lemma} \label{lem:Hext}
Every $H$-extension $X'$ of~$X$ satisfies
\begin{equation}
\MinEnt_{H'}(X') = \MinEnt_H(X)
\end{equation}
\end{lemma}

\begin{proof}
The inclusion map $i:X \hookrightarrow X'$ is $n$-monotone and $\psi = \psi' \circ i_*$. Therefore, from Lemma~\ref{lem:mon}, we have $\MinEnt_H(X) \leq \MinEnt_{H'}(X')$.

Given a metric $g$ on~$X$, we extend $g$ to the $H$-extension~$X'$ of~$X$ by induction on~$k$ so that each cell $B^k_{i_k}$ has the geometry of a long cylinder over its basis capped on its top (see~\cite[Lemma~3.5]{bab04} for a more precise construction).
Denote by~$g'$ the metric on~$X'$ so-obtained.
If the cylinders defining the geometry of the $B^k_{i_k}$'s are long enough, every loop~$\gamma'$ of~$X'$ based at some fixed point~$x_0$ of~$X$ can be homotoped in~$X'$ into a loop~$\gamma$ of~$X$ with the same based point such that $\length(\gamma) \leq \length(\gamma')$.
Since the inclusion map induces an isomorphism between $\pi_1(X)/H$ and $\pi_1(X')/H'$, we have $P_{H'}(T) \leq P_H(T)$. We deduce from Lemma~\ref{lem:entH} that $\Ent_{H'}(X',g') \leq \Ent_H(X,g)$.
Hence, $\MinEnt_{H'}(X') \leq \MinEnt_H(X)$ since the $n$-volume of the cells attached to~$X$ is null.
\end{proof}

We can now prove Theorem~\ref{theo:minent} when~$n \geq 4$. The case~$n=2$ is treated at the end of this section.\\

Let $\Phi:M \longrightarrow K$ be a degree~1 map between two oriented closed manifolds of the same dimension~$n \geq 4$ such that $K$ is a $K(\pi,1)$-space. Denote by $\phi:\pi_1(M) \longrightarrow \pi$ the homomorphism induced by $\Phi$ between the fundamental groups.
Let $\Psi:X \longrightarrow K$ be a geometric cycle representing $\Phi_*[M]=[K]$ in homology. By definition, the map~$\Psi$ is of degree~1.
Denote by $\psi:\pi_1(X) \longrightarrow~\pi$ the homomorphism induced by~$\Psi$.
The fundamental groups of~$M$ and~$X$ are related to~$\pi$ as follows.

\begin{lemma} \label{lem:surj}
The homomorphisms $\phi:\pi_1(M) \longrightarrow \pi$ and \mbox{$\psi:\pi_1(X) \longrightarrow \pi$} are epimorphisms.
\end{lemma}

\begin{proof}
Let $p:\overline{K} \longrightarrow K$ be the covering corresponding to the subgroup $\psi(\pi_1(X))$ of~$\pi_1(K) \simeq \pi$.
By definition, the fundamental group of~$\overline{K}$ is isomorphic to~$\psi(\pi_1(X))$.
The map $\Psi:X \longrightarrow K$ lifts to a map $\overline{\Psi}:X \longrightarrow \overline{K}$ uniquely defined up to deck transformations such that $\Psi=p \circ \overline{\Psi}$.
The index $[\pi:\psi(\pi_1(X))]$ of $\psi(\pi_1(X))$ in~$\pi$ is finite, otherwise $\overline{K}$ is noncompact. In this case, the homology group~$H_n(\overline{K};\Z)$ vanishes and $\Psi_*[X] = p_*(\overline{\Psi}_* [X])$ is trivial in homology, hence a contradiction.
Thus, $|\deg p| = [\pi:\psi(\pi_1(X))]$.
Since $\deg \Psi = \deg p \cdot \deg \overline{\Psi} = 1$, the index of $\psi(\pi_1(X))$ in~$\pi$ equals~$1$.
Therefore, $\psi(\pi_1(X)) = \pi$. Hence $\psi$ is surjective and the same goes for~$\phi$.
\end{proof}

The following example shows that the surjectivity of~$\psi$ does not hold anymore when $K$ is not an $n$-dimensional manifold (of course, here, we replace $[K]$ by $\Phi_*[M]$).

\begin{example}
Consider the connected sum $M = (S^1 \times S^{n-1}) \# \T^n$ with $n \geq 3$ (we can replace the $n$-torus~$\T^n$ by any essential manifold) and a $K(\pi_1(M),1)$-space $K$.
Here, $\pi_1(M) \simeq \Z * \pi_1(\T^n) \simeq \Z * \Z^n$ from Van Kampen's theorem.
Since $K$ is aspherical, the canonical map \mbox{$\Phi:M \longrightarrow K$} induced by the identity homomorphism on~$\pi_1(M)$ factorizes through $\Phi:M \longrightarrow (S^1 \times S^{n-1}) \vee \T^n \stackrel{i \vee j}{\longrightarrow} K$.
Furthermore, the map $j:\T^n \longrightarrow K$ represents~$\Phi_*[M]$ and $j_*:\pi_1(\T^n) \longrightarrow \Z * \pi_1(\T^n)$ is not surjective.
\end{example}

One can kill the kernel of~$\psi$ by attaching $2$-cells at~$X$. More precisely, we have the following.

\begin{lemma}
There exists an $H$-extension $X'$ of~$X$ and a map $\Psi':X' \longrightarrow K$ of degree~1 such that the homomorphism $\Psi'_*:\pi_1(X') \longrightarrow \pi$ induced by~$\Psi'$ is an isomorphism which agrees with~$\psi'$ defined after~\eqref{eq:X'}.

Further, since $X'$ is an $H$-extension of~$X$, $\MinEnt(X') = \MinEnt_H(X)$.
\end{lemma}

\begin{proof}
Note that the groups $\pi_1(X)$ and $\pi$ have finite presentations.
The kernel~$H$ of the epimorphism $\psi:\pi_1(X) \longrightarrow \pi$ is generated by a finite number of elements.
Define an $H$-extension~$X'$ of~$X$ by attaching $2$-cells along a finite number of (simplicial) loops representing a finite generating set of~$H$.
By construction, the map~$\Psi$ extends into a map \mbox{$\Psi':X' \longrightarrow~K$} of degree~1 such that the homomorphism induced between the fundamental groups is an isomorphism which agrees with~$\psi'$.
\end{proof}

Similar arguments lead to the following (see also~\cite[Lemme~3.11]{bab04}).

\begin{lemma}
There exists an oriented closed $n$-dimensional manifold~$M'$ and a map $\Phi':M' \longrightarrow K$ of degree~1 such that the homomorphism $\phi':\pi_1(M') \longrightarrow \pi$ induced by~$\Phi'$ is an isomorphism and $\MinEnt(M') = \MinEnt_G(M)$.
\end{lemma}

\begin{proof}
The kernel~$G$ of the epimorphism $\phi:\pi_1(M) \longrightarrow \pi$ is generated by a finite number of elements represented by disjoint simple loops $\gamma_1, \dots, \gamma_k$ of~$M$.
Consider disjoint tubular neighborhoods of the~$\gamma_i$'s homeomorphic to~$S^1 \times B^{n-1}$.
Remove these neighborhoods from~$M$ and glue $D_i^2 \times S^{n-2}$ along their boundaries~$S_i^1 \times S^{n-2}$ where $D_i^2$ is a disk bounding~$S_i^1$.
The oriented closed $n$-manifold so-obtained is noted~$M'$.

Since $K$ is aspherical and the images by~$\Phi$ of the loops~$\gamma_i$ are contractible, the map~$\Phi$ induces a map $\Phi':M' \longrightarrow K$ which represents the fundamental class of~$K$. By construction, the homomorphism $\phi':\pi_1(M') \longrightarrow K$ induced by~$\Phi'$ is an isomorphism.

Let $Y$ be the extension of~$M'$ obtained by attaching an $(n-1)$-cell $B_i^{n-1}$ along a copy $S^{n-2}$ of $D_i^2 \times S^{n-2}$ for every $1 \leq i \leq k$. This ensures us that the spheres $S^{n-2}$ of $D_i^2 \times S^{n-2}$ are contractible in~$Y$.
Therefore, there exists an $n$-monotone map $f:M \longrightarrow Y$ of degree~1 such that $\phi = \phi' \circ f_*$.
From Lemma~\ref{lem:mon} and Lemma~\ref{lem:Hext}, we have $\MinEnt_G(M) \leq \MinEnt(Y) = \MinEnt(M')$.

The reverse inequality, which will not be used in the sequel, can be obtained as follows.
Let $Z$ be the $G$-extension of~$M$ obtained by attaching $2$-cells~$D_i^2$ along the~$\gamma_i$'s. The projections of the second factor of $D_i^2 \times S^{n-2}$ onto~$D_i^2$ give rise to an $n$-monotone map of degree~1 from~$M'$ onto~$Z$ which induces an isomorphism between the fundamental groups.
As previously, from Lemma~\ref{lem:mon} and Lemma~\ref{lem:Hext}, we obtain $\MinEnt(M') \leq \MinEnt(Z) = \MinEnt_G(M)$.
\end{proof}

The following result rests on a construction of~\cite[Th\'eor\`eme~3.7]{bab04}.
The restriction on the dimension, i.e., $n \neq 3$, is required at this point precisely (see Remark~\ref{rem}).

\begin{lemma} \label{lem:key}
There exists an extension $X''$ of $X'$ and a degree~1 map $\Phi'':M' \longrightarrow X''$ such that the homomorphism \mbox{$\phi'':\pi_1(M') \longrightarrow \pi_1(X'')$} induced by~$\Phi''$ is an isomorphism which agrees with~$\psi'^{-1} \circ \phi'$ and 
\begin{equation}
\MinEnt(X'') \geq \MinEnt(M').
\end{equation}

Further, since $X''$ is an extension of~$X'$, $\MinEnt(X'') = \MinEnt(X')$.
\end{lemma}

\begin{proof}
By assumption, $\psi':\pi_1(X') \longrightarrow \pi$ and $\phi':\pi_1(M') \longrightarrow \pi$ are isomorphisms and $\Psi'_*[X']=\Phi'_*[M']=[K]$ in $H_n(K;\Z) \simeq H_n(\pi;\Z)$.
We can apply the theorem~3.7 of~\cite{bab04} which yields an extension~$X''$ of~$X'$ and a map $\Phi'':M' \longrightarrow X''$ such that the homomorphism \mbox{$\phi'':\pi_1(M') \longrightarrow \pi_1(X'')$} induced by~$\Phi''$ satisfies~$\phi'' = {\psi'}^{-1} \circ \phi'$.

\begin{remark} \label{rem}
Strictly speaking, this result has been stated for two manifolds $M_1$ and $M_2$ (see~\cite[Th\'eor\`eme~3.7]{bab04} for the notations).
Furthermore, it does not hold anymore when $M_2$ is replaced by a polyhedron (see~\cite{bab04} for a counterexample).
However, when $M_1$ is a pseudomanifold and $M_2$ is a manifold, the proof still applies which permits us to conclude.

Let us develop this point by recalling the steps of the construction of~\cite[Th\'eor\`eme~3.7]{bab04}:
\begin{enumerate}
\item[1)] The homotopy groups of~$M_{1}$ of dimension~$i$ with $2 \leq i \leq n-1$ can be successively killed by attaching cells of dimensions~$i+1$.
The resulting spaces are CW-complexes noted~$M_{1}(i+1)$ with $M_{1}(n)$ homotopically equivalent to the $n$-skeleton~$K(\pi,1)^{(n)}$ of~$K(\pi,1)$.
\item[2)] The cellular approximation theorem applied to the classifying map of~$M_{2}$ yields a map $g:M_{2} \longrightarrow M_{1}(n)$.
Note that $M_{1}(n)$ is not an extension of~$M_{1}$.
\item[3)] The map~$g$ is not uniquely defined up to homotopy.
However, it may be chosen so that it can be deformed into a map $\hat{g}:M_{2} \longrightarrow M_{1}(n-1)$ inducing the same homomorphism as~$g$ between the fundamental groups.
\item[4)] By compactness, the map~$\hat{g}$ can be deformed into a map whose image is a finite subcomplex of~$M_{1}(n-1)$ and so an extension of~$M_{1}$.
\end{enumerate}

Step~3) is the main step of the construction.
It relies on the obstruction to extend a map inducing a fundamental group isomorphism.
When such an obstruction is trivial, the extension can be performed without modifying the original map on its $2$-codimensional skeleton and, in particular, on its $2$-skeleton when~$n \geq 4$.
This ensures that the extension still induces an isomorphism between the fundamental groups.

The proof of~3) requires that $M_{2}$ is a manifold to apply the Poincar\'e duality.
However, the assumption that $M_{1}$ is a manifold plays no role in any of the four steps of the construction ($M_{1}$ is turned into a CW-complex from the first step by attaching cells to it).
For the purpose of the proof, it is enough that $M_{1}$ is a pseudomanifold.
\end{remark}

Let us now conclude the proof of Lemma~\ref{lem:key}.
Since $K$ is aspherical, the map~$\Psi'$ extends to~$X''$ into a degree~1 map.
This map, still denoted~$\Psi'$, induces a homomorphism between the fundamental groups which agrees with $\psi':\pi_1(X'') \simeq \pi_1(X') \longrightarrow \pi$.
The homomorphism~$\psi' \circ \phi''$ agrees with~$\phi'$. Therefore, since $K$ is aspherical, the map $\Psi' \circ \Phi''$ is homotopic to the degree~1 map~$\Phi'$. This shows that $\deg \Psi' \circ \Phi'' = \deg \Psi' \cdot \deg \Phi'' = 1$.
Thus, the degree of~$\Phi''$ is equal to~1.

The polyhedron~$X''$ is an $H$-extension of the closed $n$-pseudomanifold~$X$. Therefore, it admits a cellular decomposition with exactly one $n$-cell (see Lemma~\ref{lem:cell}).
As noticed in~\cite[Propri\'et\'e~3.12]{bab04}, the proof of~\cite[Theorem~4.1]{eps} still applies in this case. This implies from~\cite[Propri\'et\'e~3.12]{bab04} that the map $\Phi'':M' \longrightarrow X''$ of degree~1 is homotopic to an \mbox{$n$-monotone} map $\overline{\Phi}'':M' \longrightarrow~X''$.
We deduce from Lemma~\ref{lem:mon} that $\MinEnt(X'') \geq \MinEnt(M')$.
\end{proof}

The three previous lemmas show that every geometric cycle~$X$ representing~$\Phi_*[M]$ satisfies
\begin{equation}
\MinEnt_G(M) \leq \MinEnt_H(X).
\end{equation}
Thus, $\MinEnt_G(M) \leq \MinEnt(\Phi_*[M])$.

The reverse inequality is obvious. Therefore, Theorem~\ref{theo:minent} is proved for~$n \geq 4$. \\

When $n=2$, $M$ is a surface of genus~$\gamma \geq 1$. In this case, $K$ is homeomorphic to a surface of genus at most~$\gamma$. Let $\Psi:X \longrightarrow K$ be a geometric cycle representing~$\Phi_*[M]=[K]$ in homology.
From the definition of pseudomanifolds, there exists a finite number of points $x_1, \dots, x_k$ of~$X$ such that $X \setminus \{ x_1, \dots, x_k \}$ is a two-dimensional manifold and every point~$x_i$ admits a neighborhood homeomorphic to a cone over a finite number of circles. The completed space of $X \setminus \{ x_1, \dots, x_k \}$ is a closed orientable surface~$N$.
The inclusion $X \setminus \{ x_1, \dots, x_k \} \hookrightarrow X$ extends into a $2$-monotone map $j:N \longrightarrow X$ such that $j_*[N]=[X]$.
Thus, the map $\Psi_0=\Psi \circ j :N \longrightarrow K$ represents~$\Phi_*[M]=[K]$.
We can assume that the genus of~$K$ is greater than one, otherwise the result is obvious since $\pi_1(M)/G$ is isomorphic to~$\Z^2$ and the minimal $G$-entropy of~$M$ is then zero.
The degree~1 map~$\Psi_0$ induces a homomorphism $\psi_0:\pi_1(N) \longrightarrow \pi$ between the fundamental groups, whose kernel is noted~$H_0$.
Arguing as in~\cite[Th\'eor\`eme~8.1]{bcg95}, one can show that $\MinEnt_{H_0}(N) = \MinEnt(M)$.
Since $j$ is $2$-monotone and $\psi_0 = \psi \circ j_*$, we deduce that $\MinEnt_{H_0}(N) \leq \MinEnt_H(X)$ from Lemma~\ref{lem:mon}.
Therefore, $\MinEnt(M) \leq \MinEnt_H(X)$. As previously, Theorem~\ref{theo:minent} immediately follows when~$n=2$.

\section{Entropy of extremal geometric cycles} \label{sec:ent}

We introduce in this section the notion of regular geometric cycles and bound from above their relative entropy in terms of their relative systolic volume. \\

In dimension greater than two, we do not know whether the infimum in inequality~\eqref{eq:sigma} is achieved by any metric. The likely absence of extremal metrics in general makes the study of systolic geometry delicate.
In order to get around this difficulty M.~Gromov proved that every essential manifold can be represented by a regular geometric cycle (see below). These regular geometric cycles share enough patterns we expect from possible extremal metrics to play the role of substitute in their absence.
Note however that regular geometric cycles are not necessarily homeomorphic to the manifold they represent, neither do they necessarily have the same fundamental group. \\

Given a discrete group~$\pi$, the optimal systolic volume of a homology class $h \in H_n(\pi;\Z)$ is defined as
\begin{equation}
\sigma(h) = \inf_X \sigma_\psi(X)
\end{equation}
where $X=(X,\Psi,g)$ runs over all the geometric cycles representing~$h$.

M.~Gromov showed in~\cite[\S6]{gro83} that there exists a positive constant~$c_n$ such that every nontrivial homology class $h \in H_n(\pi;\Z)$ satisfies
\begin{equation}
\sigma(h) \geq c_n
\end{equation}
Furthermore, he proved in~\cite[p.~71]{gro83} that there exists $c_n > 0$ such that, for every $\varepsilon >0$, every nontrivial homology class $h \in H_n(\pi;\Z)$ can be represented by a geometric cycle~$X=(X,\Psi,g)$ such that
\begin{equation} \label{eq:reg0}
\sigma_\psi(X,g) \leq (1+\varepsilon) \sigma(h)
\end{equation}
and
\begin{equation} \label{eq:reg}
\vol_g(B(x,R)) \geq c_n R^n
\end{equation}
for every $x \in X$ and $\varepsilon \leq R / \sys_\psi(X,g) \leq \frac{1}{2}$, where $\vol_g(B(x,R))$ is the volume of the $g$-ball centered at~$x$ with radius~$R$ in~$X$. In~\eqref{eq:reg}, we can take $c_n = n^{-n+2} ((n-1)! \sqrt{n!})^{-n+1}$.

Such geometric cycles are said to be $\varepsilon$-regular. \\

Arguing as in~\cite{k-s} (see also~\cite{kat}), we show the following.

\begin{proposition} \label{prop:Hent}
Every $\varepsilon$-regular geometric cycle $(X,\Psi,g)$ representing a nontrivial homology class $h \in H_n(\pi;\Z)$ satisfies
\begin{equation} 
\label{eq:Hent}
\Ent_H(X,g) \leq \frac{1}{\beta \sys_\psi(X,g)} \log \left( \frac{\sigma_\psi(X,g)}{c_n \alpha^n} \right),
\end{equation}
whenever $\alpha \geq \varepsilon$, $\beta >0$ and $4\alpha + \beta < \frac{1}{2}$, where~$c_n$ is given by~\eqref{eq:reg}.
\end{proposition}

\begin{proof}
Let $x_0\in X$ be a fixed basepoint.  Consider a maximal system of
disjoint balls
\begin{equation}
\label{22}
B_i= B(x_i,R) \subset X
\end{equation}
of radius $R=\alpha \sys_\psi(X,g)$ and centers~$x_i$ with~$i \in I$,
including $x_0$.  Since the geometric cycle $(X,\Psi,g)$ is $\varepsilon$-regular, we have
\begin{equation} \label{eq:*}
\vol B_i \geq c_n \alpha^n \sys_\psi(X,g)^n \quad \forall i\in I.
\end{equation}
Therefore, this system admits at most $\frac{\vol(X,g)}{c_n \alpha^n
\sys_\psi(X,g)^n}$ balls. Thus,
\begin{equation} \label{eq:**}
|I| \leq \frac{\sigma_\psi(X,g)}{c_n \alpha^n}.
\end{equation}
Let $c:[0,T] \to X$ be a geodesic loop of length~$T$ based at~$x_0$ whose image by~$\Psi$ is noncontractible.
Let
\begin{equation}
\label{23}
m = \left[\frac{T}{\beta \sys_\psi(X,g)} \right]
\end{equation}
be the integer part.  The point $p_0=x_0$, together with the points
\[
p_k = c(k \beta \sys_\psi(X,g)), \quad k=1, \ldots, m
\]
and the point $p_{m+1}=x_0$, partition the loop~$c$ into $m+1$
segments of length at most~$\beta \sys_\psi(X,g)$.  Since the system of
balls~$B_i$ is maximal, the disks of radius $2R = 2 \alpha
\sys_\psi(X,g)$ centered at~$x_i$ cover~$X$.  Therefore, for every~$p_k$,
a nearest point~$q_k$ among the centers~$x_i$, is at distance at
most~$2R$ from~$p_k$.  Consider the loop
\[
\alpha_k= c_k \cup [p_{k+1},q_{k+1}] \cup [q_{k+1},q_k] \cup
[q_k,p_k],
\]
where $c_k$ is the arc of~$c$ joining $p_k$ to~$p_{k+1}$, while
$[x,y]$ denotes a minimizing path joining $x$ to~$y$.  Then
\[
\length(\alpha_k) \leq 2(4 \alpha + \beta) \sys_\psi(X,g) < \sys_\psi(X,g),
\]
by our hypothesis on $\alpha, \beta$.  Thus the image by~$\Psi$ of the loops $\alpha_k$ is contractible.  The same is true for the loops $c_0 \cup [p_1,q_1] \cup
[q_1,x_0]$ and $c_m \cup [x_m,q_m] \cup [q_m,p_m]$.

Therefore, the image by~$\Psi$ of the geodesic loop~$c$ is homotopic to the image by~$\Psi$ of a piecewise geodesic loop~$c'= (x_0,q_1, \dots, q_{m} ,x_0)$.  Note that since the length of minimizing paths from $p_k$ to $q_k$ is at most $(4 \alpha + \beta) \sys_\psi(X,g) < \frac{1}{2} \sys_\psi(X,g)$, their images by~$\Psi$ are homotopic with fixed endpoints. Thus, the choice of a minimizing path between~$p_k$ and~$q_k$ does not matter.

Thus, $c$ and $c'$ represent the same class in~$\Aut(X_H) \simeq \pi_1(X)/H$ where $H= \ker \psi$.
Furthermore, two closed geodesics $c_1$ and $c_2$ which induce two distinct classes in~$\Aut(X_H)$ give rise to two loops $c'_1$ and $c'_2$ with the same property.
Thus, the number $P_H(T)$ of classes in~$\Aut(X_H)$ which can be represented by loops of length at most~$T$ based at~$x_0$ satisfies
\begin{eqnarray}
P_H(T) & \leq & |I|^{m} \nonumber \\ & \leq & |I|^{\frac{T}{\beta
      \sys_\psi(X,g)}} \nonumber \\ & \leq & \left( \frac{\sigma_\psi(X,g)}{c_n \alpha^n} \right)^{\frac{T}{\beta
      \sys_\psi(X,g)}}, \label{eq:***}
\end{eqnarray}
and the proposition now follows from Lemma~\ref{lem:entH}.
\end{proof}

\section{Relative systolic volume and relative minimal entropy} \label{sec:Gsys}

Using the results established in the previous sections, we show that under some topological conditions, the relative systolic volume of manifolds is bounded from below in terms of their relative minimal entropy. More precisely, we have the following.

\begin{theorem} \label{theo:Gsys}
Let $\pi$ be a discrete group.
Let $\Phi:M \longrightarrow K$ be a degree~1 map between two oriented closed manifolds of the same dimension~$n \neq 3$ such that $K$ is a $K(\pi,1)$-space.
Then, there exists a positive constant~$C_n$ depending only on~$n$ such that
\begin{equation} \label{eq:Gsys}
\sigma_\phi(M) \geq C_n \frac{\MinEnt_G(M)^n}{\log^n(1+\MinEnt_G(M))}.
\end{equation}
where $G$ is the kernel of the homomorphism $\phi:\pi_1(M) \longrightarrow \pi$ induced by~$\Phi$.
\end{theorem}

\begin{remark}
Since $\sigma(M) \geq \sigma_\phi(M)$, we can replace $\sigma_\phi(M)$ by $\sigma(M)$ in inequality~\eqref{eq:Gsys}.
\end{remark}

\begin{example}
Recall M.~Gromov's isosystolic inequality for surfaces of large genus (see~\cite[6.4.D']{gro83}, \cite{k-s} and \cite{bal}): there exists a positive constant~$c$ such that every closed surface~$\Sigma_{\gamma}$ of genus~$\gamma>1$ satisfies
$$\sigma(\Sigma_{\gamma}) \geq c \frac{\gamma}{(\log \gamma)^2}.$$

Theorem~\ref{theo:Gsys} generalizes this result for surfaces of large genus with few small handles. More precisely, consider a degree~1 map $\Phi: \Sigma_{\gamma} \longrightarrow \Sigma_{\gamma'}$ between two closed surfaces of genus $\gamma$ and $\gamma'>1$.
From Remark~\ref{rem:minent} and inequality~\eqref{eq:bcg}, we have
$$\MinEnt_G(\Sigma_{\gamma}) = \MinEnt(\Sigma_{\gamma'}) = 4 \pi(\gamma' -1).$$
This implies that 
$$\sigma_\phi(\Sigma_{\gamma}) \geq c' \frac{\gamma'}{(\log \gamma')^2}$$
for some~$c' >0$.
Thus, the area of a closed Riemannian surface of large genus with, let's say, a single ''small'' handle is still large.
\end{example}

\begin{proof}[Proof of Theorem~\ref{theo:Gsys}]
We will use the same notations as in the previous sections.
From Proposition~\ref{prop:Hent}, every $\varepsilon$-regular cycle $(X,\Psi,g)$ representing the image~$\Phi_*[M]$ of the fundamental class of~$M$ in $H_n(K;\Z) \simeq H_n(\pi;\Z)$ satisfies
\begin{equation}
\Ent_H(X,g)^n \vol(X,g) \leq \frac{\sigma_\psi(X,g)}{\beta^n} \log^n \left( \frac{\sigma_\psi(X,g)}{c_n \alpha^n} \right)
\end{equation}
whenever $\alpha \geq \varepsilon$, $\beta > 0$ and $4 \alpha + \beta < \frac{1}{2}$, with $\sigma_\psi(X,g) \leq (1+\varepsilon) \sigma_\phi(M)$.
From Theorem~\ref{theo:minent}, we have
\begin{equation}
\MinEnt_G(M)^n = \MinEnt(\Phi_*[M])^n \leq \Ent_H(X,g)^n \vol(X,g).
\end{equation}
Therefore, we deduce that
\begin{equation} \label{eq:minent}
\MinEnt_G(M)^n \leq \frac{\sigma_\phi(M)}{\beta^n} \log^n \left( \frac{\sigma_\phi(M)}{c_n \alpha^n} \right)
\end{equation}
whenever $\alpha, \beta > 0$ and $4 \alpha + \beta \leq \frac{1}{2}$.

Define
\[
\alpha =
\begin{cases}
\frac{1}{10} & \text{if } n^n c_n \leq 1 \\
\frac{1}{10 n \sqrt[n]{c_n}} & \text{otherwise} \\
\end{cases}
\text{ and }
\beta =
\begin{cases}
\frac{1}{10} n \sqrt[n]{c_n} & \text{if } n^n c_n \leq 1 \\
\frac{1}{10} & \text{otherwise.} \\
\end{cases}
\]
In both cases, we have $4 \alpha + \beta \leq \frac{1}{2}$ and $\beta = n \alpha \sqrt[n]{c_n}$.

Let $\rho = \left( \frac{\sigma_\phi(M)}{c_n \alpha^n} \right)^{\frac{1}{n}}$ and $\delta = \frac{\beta}{n \alpha \sqrt[n]{c_n}} \MinEnt_G(M)$.
For our choices of $\alpha$ and $\beta$, we have $\delta = \MinEnt_G(M)$. Note also that~$\rho \geq 1$ since $\sigma_\phi(M) \geq \frac{c_n}{2^n}$ from inequality~\eqref{eq:reg}.

For $\delta \geq e$, inequality~\eqref{eq:minent} yields the following estimate:
\[
\begin{aligned}
\rho \log \rho & \geq \delta \\
               & \geq \delta - \frac{\delta \log \log \delta}{\log \delta} \\
               & = \frac{\delta}{\log \delta} \log \left( \frac{\delta}{\log \delta} \right).
\end{aligned}
\]
The function $x \log x$ is increasing for $x \geq \frac{1}{e}$.
Thus, since $\rho \geq 1$ and $\frac{\delta}{\log \delta} \geq e$ for $\delta \geq e$, we deduce that $\rho \geq \frac{\delta}{\log\delta} \geq \frac{\delta}{\log(1+\delta)}$ for $\delta \geq e$.

For $\delta \leq e$, we have $0 \leq \frac{\delta}{\log(1+\delta)} \leq \frac{e}{\log(1+e)}$. Therefore, we obtain
$\rho \geq 1 \geq \frac{\log(1+e)}{e} \frac{\delta}{\log(1+\delta)}$ for $\delta \leq e$.

In conclusion, for every $\delta \geq 0$, the following holds:
\begin{equation}
\rho \geq \frac{\log(1+e)}{e} \frac{\delta}{\log(1+\delta)}.
\end{equation}
Hence,
\begin{equation}
\sigma_\phi(M) \geq c_n \alpha^n \frac{\log^n(1+e)}{e^n} \frac{\MinEnt_G(M)^n}{\log^n(1+\MinEnt_G(M))}.
\end{equation}
Therefore, inequality~\eqref{eq:Gsys} holds with~$C_n = c_n \alpha^n \frac{\log^n(1+e)}{e^n}$.
\end{proof}

\section{Relative systolic volume and algebraic entropy} \label{sec:Galg}

In this section, we show how the relative entropy of a polyhedron is related to the algebraic entropy of a quotient of its fundamental group.
Then, we derive a lower bound on the relative systolic volume of some manifolds in terms of the algebraic entropy of some groups. Specifically, we prove the following.

\begin{theorem} \label{theo:Galg}
Let $\pi$ be a discrete group.
Let $\Phi:M \longrightarrow K$ be a degree~1 map between two oriented closed manifolds of the same dimension~$n$ such that $K$ is a $K(\pi,1)$-space.
Then, there exists a positive constant~$C_n$ depending only on~$n$ such that
\begin{equation} \label{eq:Galg}
\sigma_\phi(M) \geq C_n \frac{\Ent_{alg}(\pi)}{\log(1+\Ent_{alg}(\pi))}
\end{equation}
where $\phi:\pi_{1}(M) \longrightarrow \pi$ is the homomorphism induced by~$\Phi$.
\end{theorem}

\begin{remark}
We can replace $\sigma_\phi(M)$ by $\sigma(M)$ in inequality~\eqref{eq:Galg}.
Furthermore, the integer~$n$ can be equal to~3.
\end{remark}

\begin{remark}
Let $F_k$ be the free group with $k$ generators, $\pi_1(\Sigma_{\gamma})$ be the fundamental group of a closed surface of genus~$\gamma$ and $\pi$ be a discrete group with $k$ generators and $p$ relations. Then, we have
\begin{itemize}
\item $\Ent_{alg}(F_k) \geq \log(2k-1)$ (see~\cite[5.13]{gro99} or~\cite[VII.B]{har});
\item $\Ent_{alg}(\pi_1(\Sigma_{\gamma})) \geq \log(4\gamma-3)$ (see~\cite[VII.B]{har});
\item $\Ent_{alg}(\pi) \geq c \log(2(k-p)-1)$ for some~$c>0$ (see~\cite[5.14]{gro99} and~\cite[p.~82-83]{gro81}).
\end{itemize}
\end{remark}

As shown in the following example, Theorem~\ref{theo:Galg} may yield nontrivial lower bounds on the systolic volume when the contents of Theorem~\ref{theo:simpl} and Theorem~\ref{theo:ent} are trivial.

\begin{example}
Consider the closed aspherical manifolds
$$M_{\gamma} = \Sigma_{\gamma} \times \T^{n-2}$$
of dimension~$n \geq 3$, where $\Sigma_{\gamma}$ is the closed surface of genus~$\gamma$ and $\T^{n-2}$ is the $(n-2)$-torus. The group~$\pi_1(\Sigma_{\gamma})$ is a quotient of~$\pi_1(M)$. Therefore,
$$\Ent_{alg}(\pi_1(M_{\gamma})) \geq \Ent_{alg}(\Sigma_{\gamma}) \geq \log(4\gamma-3).$$
From Theorem~\ref{theo:Galg}, we obtain
$$\sigma(M_{\gamma}) \gtrsim \frac{\log \gamma}{\log \log \gamma},$$
which goes to infinity as $\gamma \rightarrow \infty$.
Such an estimate follows neither from Theorem~\ref{theo:simpl} nor from Theorem~\ref{theo:ent} since $\MinEnt(M_{\gamma}) = ||M_{\gamma}|| = 0$.
\end{example}

In this second example, we show that the conclusion of Theorem~\ref{theo:Galg} does not hold without some topological conditions on~$M$.

\begin{example} \label{ex:handles}
Fix a closed aspherical manifold~$M$ of dimension~\mbox{$n \geq 3$}.
Consider the sequence of essential manifolds
$$M_k = M \# S^1 \times S^{n-1} \# \dots \# S^1 \times S^{n-1}$$
where we take the connected sum of~$M$ with $k$ copies of~$S^1 \times S^{n-1}$.
As shown in~\cite{b-b}, adding handles to a manifold does change the value of the optimal systolic volume. Thus, for every~$k$, we have
$$\sigma(M_k) = \sigma(M).$$
On the other hand, the difference between the number of generators and the number of relations of~$\pi_1(M_k)$ is equivalent to~$k$. Therefore,
$$\Ent_{alg}(\pi_1(M_k)) \gtrsim \log k.$$
This shows that inequality~\eqref{eq:Galg} does not hold in general: some topological condition must be satisfied as in Theorem~\ref{theo:Galg}.
\end{example}

The following result is a relative version of~\cite[Theorem~5.16]{gro99} (which is incorrectly stated).

\begin{lemma} \label{lem:alg}
Let $X$ be a closed polyhedron and $H$ be a normal subgroup of the fundamental group~$\pi_1(X)$ of~$X$. Then, every piecewise linear Riemannian metric~$g$ on~$X$ satisfies
\begin{equation}
\Ent_{alg}(\pi_1(X)/H) \leq 2 \Diam(X,g) \Ent_H(X,g).
\end{equation}
\end{lemma}

\begin{proof}
Fix~$x_0$ in~$X$. The fundamental group~$\pi_1(X,x_0)$ of~$X$ is generated by a finite set~$\Sigma$ formed of classes represented by loops based at~$x_0$ of length at most $2 \Diam(X,g)$ (see~\cite[3.22]{gro99}). Denote by~$\overline{\Sigma}$ the generating set of $\Gamma:=\pi_1(X,x_0)/H$ induced by~$\Sigma$.
Every class of $|.|_{\overline{\Sigma}}$-length at most~$k$ can be represented by a loop based at~$x_0$ of length at most~$2k \Diam(X,g)$.
Therefore, $N_{\overline{\Sigma}}(k) \leq P_H(2k\Diam(X,g))$ where $N_{\overline{\Sigma}}$ and $P_H$ are defined in~\eqref{eq:entalg} and Lemma~\ref{lem:entH}.
Thus,
\begin{equation}
\frac{\log(N_{\overline{\Sigma}}(k))}{k} \leq 2 \Diam(X,g) \frac{\log(P_H(2k\Diam(X,g)))}{2k\Diam(X,g)}.
\end{equation}
From Lemma~\ref{lem:entH}, we obtain
\begin{equation}
\Ent_{alg}(\Gamma) \leq \Ent_{alg}(\Gamma,\overline{\Sigma}) \leq 2 \Diam(X,g) \Ent_H(X,g).
\end{equation}
when $k$ goes to infinity.
\end{proof}

We can now prove the main result of this section.

\begin{proof}[Proof of Theorem~\ref{theo:Galg}]
Every $\varepsilon$-regular geometric cycle~$(X,\Psi,g)$ representing the image~$\Phi_*[M]$ of the fundamental class of~$M$ admits at least $\left[ \frac{\Diam(X,g)}{\sys_\psi(X,g)} \right]$ disjoint balls of radius~$\frac{1}{2} \sys_\psi(X,g)$.
Combined with the inequality~\eqref{eq:reg}, this yields
\begin{equation}
\vol(X,g) \geq \left[ \frac{\Diam(X,g)}{\sys_\psi(X,g)} \right] \frac{c_n}{2^n} \sys_\psi(X,g)^n.
\end{equation}
Hence,
\begin{equation} \label{eq:6*}
\left[ \frac{\Diam(X,g)}{\sys_\psi(X,g)} \right] \leq \frac{2^n}{c_n} \sigma_\psi(X,g)
\end{equation}
From Lemma~\ref{lem:surj} (which holds for every integer~$n$), we have $\pi_1(X)/H \simeq~\pi$.
Now, Lemma~\ref{lem:alg} and Proposition~\ref{prop:Hent} imply that
\begin{eqnarray}
\Ent_{alg}(\pi) & \leq & 2 \Diam(X,g) \Ent_H(X,g) \\
             & \leq & \frac{2 \Diam(X,g)}{\beta \sys_\psi(X,g)} \log \left( \frac{\sigma_\psi(X,g)}{c_n \alpha^n} \right) \label{eq:6**}
\end{eqnarray}
whenever $\alpha,\beta >0$ and $4\alpha + \beta < \frac{1}{2}$.

If $\frac{\Diam(X,g)}{\sys_\psi(X,g)} \geq 2$, then $\frac{\Diam(X,g)}{\sys_\psi(X,g)} \leq 2 \left[ \frac{\Diam(X,g)}{\sys_\psi(X,g)} \right]$. \\
Otherwise, $\frac{\Diam(X,g)}{\sys_\psi(X,g)} \leq \frac{2^{n+1}}{c_n} \sigma_\psi(X,g)$ since $\sigma_\psi(X,g) \geq \frac{c_n}{2^n}$ from~\eqref{eq:reg}.

In both case, we deduce from~\eqref{eq:6*} and~\eqref{eq:6**} that
\begin{equation}
\Ent_{alg}(\pi) \leq \frac{2^{n+1}}{\beta c_n} \sigma_\psi(X,g) \log \left( \frac{\sigma_\psi(X,g)}{c_n \alpha^n} \right)
\end{equation}
whenever $\alpha,\beta >0$ and $4\alpha + \beta < \frac{1}{2}$, with $\sigma_\psi(X,g) \leq (1+\varepsilon) \sigma_\phi(M)$.

Therefore, we have
\begin{equation}
\Ent_{alg}(\pi) \leq \frac{2^{n+1}}{\beta c_n} \sigma_\phi(M) \log \left( \frac{\sigma_\phi(M)}{c_n \alpha^n} \right)
\end{equation}
whenever $\alpha,\beta >0$ and $4\alpha + \beta \leq \frac{1}{2}$.

Arguing as in the proof of Theorem~\ref{theo:Gsys}, we obtain
\begin{equation}
\sigma_\phi(M) \geq c_n \alpha^n \frac{\log(1+e)}{e} \frac{\Ent_{alg}(\pi)}{\log(1+\Ent_{alg}(\pi))}
\end{equation}
with $\alpha=\frac{1}{10}$. Therefore, inequality~\eqref{eq:Galg} holds with $C_n = c_n \alpha^n \frac{\log(1+e)}{e}$.
\end{proof}

\end{document}